\documentclass[12pt, reqno, twoside]{amsart}
\usepackage{amsmath, amsthm, mathrsfs, amscd, amsfonts, amssymb, graphicx, color}
\usepackage[bookmarksnumbered, colorlinks, plainpages]{hyperref}
\hypersetup{colorlinks=true,linkcolor=red, anchorcolor=green, citecolor=cyan, urlcolor=red, filecolor=magenta, pdftoolbar=true}
\usepackage{enumitem}
\usepackage{enumitem}
\newtheorem{theorem}{Theorem}[section]
\newtheorem{lemma}[theorem]{Lemma}

\theoremstyle{definition}

\theoremstyle{remark}
\newtheorem{remark}[theorem]{Remark}
\numberwithin{equation}{section}

\begin{document}

\setcounter{page}{1}

\title[Gradient estimates for nonlinear elliptic equations]{Gradient estimates for nonlinear elliptic equations involving 
the Witten Laplacian on smooth metric measure spaces and implications}

\author[A. Taheri]{Ali Taheri}

\author[V. Vahidifar]{Vahideh Vahidifar}

\address{School of Mathematical and  Physical Sciences, 
University of Sussex, Falmer, Brighton, United Kingdom.}
\email{\textcolor[rgb]{0.00,0.00,0.84}{a.taheri@sussex.ac.uk}} 
\address{School of Mathematical and  Physical Sciences, 
University of Sussex, Falmer, Brighton, United Kingdom.}
\email{\textcolor[rgb]{0.00,0.00,0.84}{v.vahidifar@sussex.ac.uk}}


\begin{abstract}
This article presents new local and global gradient estimates of Li-Yau type for positive solutions to a class of nonlinear elliptic 
equations on smooth metric measure spaces involving the Witten Laplacian. The estimates are derived 
under natural lower bounds on the associated Bakry-\'Emery Ricci curvature tensor and find utility in 
proving general Harnack inequalities and Liouville-type theorems to mention a few. 
The results here unify, extend and improve various existing results in the literature for special nonlinearities 
already of huge interest and applications. Some important consequences are presented and discussed. 
\end{abstract}

\maketitle
{ 
\hypersetup{linkcolor=black}
\tableofcontents
}

\section{Introduction} \label{sec1}

Suppose $(M, g, d\mu)$ is a smooth metric measure space by which it is meant that $M$ is a complete 
Riemannian manifold of dimension $n \ge 2$ endowed with a weighted measure $d\mu=e^{-f} dv_g$, 
$f$ is a smooth potential on $M$, $g$ is the Riemannian metric and $dv_g$ is the usual 
Riemannain volume measure. In this paper we derive 
gradient estimates of local and global Li-Yau type along with Harnack inequalities and Liouville-type results for 
positive smooth solutions $u$ to the nonlinear elliptic equation: 
\begin{align} \label{eq11}
\Delta_f u (x) + \Sigma [x, u(x)] =0,  \qquad \Sigma: M \times \mathbb{R} \to \mathbb{R}. 
\end{align}

Here $\Delta_f u = e^f\text{div}[e^{-f} \nabla u] = \Delta u - \langle \nabla f, \nabla u \rangle$ is the Witten 
Laplacian (also known as the weighted or drifting Laplacian, or occasionally to emphasise the choice of 
$f$, the $f$-Laplacian) where $\nabla, {\rm div}$ and $\Delta$ are the usual gradient, 
divergence and Laplace-Beltrami operators associated with the metric $g$. 
The Witten Laplacian is a symmetric diffusion operator with respect to the invariant weighted measure 
$d\mu=e^{-f} dv_g$ and arises in many contexts ranging from probability theory, geometry 
and stochastic processes to quantum field theory and statistical mechanics \cite{BD, Bak, Gr}. 
It is the natural generalisation of the Laplace-Beltrami operator to the smooth metric measure 
space setting and it coincides with the latter precisely when the potential $f$ is a constant.

The nonlinearity $\Sigma=\Sigma[x,u]$ in \eqref{eq11} 
is a sufficiently smooth multi-variable function depending on both the spatial variable $x \in M$ and the independent 
variable (solution) $u$. In what follows, in order to better orient the reader and showcase the results, we discuss 
various examples of such nonlinearities arising from different contexts, e.g., from conformal geometry and 
mathematical physics, each presenting a different phenomenon whilst depicting a corresponding singular 
or regular behaviour.

As for the curvature properties of the triple $(M, g, d\mu)$ we introduce the generalised Bakry-\'Emery Ricci curvature tensor 
({\it see} \cite{BD, BE, Lott}), by writing, for $m \ge n$,  
\begin{equation} \label{Ricci-m-f-intro}
{\mathscr Ric}^m_f(g) = {\mathscr Ric}(g) + {\rm Hess}(f) - \frac{\nabla f \otimes \nabla f}{m-n}. 
\end{equation}

Here ${\mathscr Ric}(g)$ denotes the Riemannain Ricci curvature tensor of $g$, ${\rm Hess}(f)$ stands for the 
Hessian of $f$, and $m \ge n$ is a fixed constant. For the sake of clarity we point out that in the event 
$m=n$, by convention, $f$ is only allowed to be a constant, resulting in ${\mathscr Ric}^n_f(g)={\mathscr Ric}(g)$, 
whilst, we also allow for $m = \infty$ in which case by formally passing to the limit in \eqref{Ricci-m-f-intro} we write 
${\mathscr Ric}_f^\infty(g) = {\mathscr Ric}(g) + {\rm Hess}(f) := {\mathscr Ric}_f(g)$.

According to the weighted Bochner-Weitzenb\"ock formula, for every smooth function $u$ on $M$ we have, 
\begin{equation} \label{Bochner}
\frac{1}{2} \Delta_f |\nabla u|^2 = |{\rm Hess}(u)|^2 + \langle \nabla u, \nabla \Delta_f u \rangle + {\mathscr Ric}_f (\nabla u, \nabla u).
\end{equation}

Hence by an application of the Cauchy-Schwartz inequality giving $\Delta u \le \sqrt n |{\rm Hess}(u)|$ and upon recalling the identity 
$\Delta_f u = \Delta u - \langle \nabla f, \nabla u \rangle$ it is evident that 
\begin {equation} \label{CS-ineq-Eq1.4}
|{\rm Hess}(u)|^2 \ge \frac{(\Delta u)^2}{n}, \qquad  
 \frac{(\Delta u)^2}{n} + \frac{\langle \nabla f, \nabla u \rangle ^2}{m-n} \ge \frac{(\Delta_f u)^2}{m},  
\end{equation}
and so it follows from \eqref{Bochner} and \eqref{CS-ineq-Eq1.4} that 
\begin{equation}
\frac{1}{2} \Delta_f |\nabla u|^2 - \langle \nabla u, \nabla \Delta_f u \rangle 
\ge \frac{1}{m} (\Delta_f u)^2  + {\mathscr Ric}^m_f (\nabla u, \nabla u).
\end{equation}
In particular, subject to a curvature lower bound ${\mathscr Ric}_f^m(g) \ge \mathsf{k} g$, the operator $L=\Delta_f$ is seen to satisfy the 
curvature-dimension condition ${\rm CD}(\mathsf{k},m)$ ({\it cf.} \cite{BD, BE, Bak, VC}).

Our principal objective in this paper is to develop local and global gradient estimates of Li-Yau type and Harnack inequalities for positive smooth 
solutions to \eqref{eq11}. It is well-known that these estimates form the basis for deriving various 
qualitative properties of solutions and are thus of great significance and utility ({\it see}, e.g., \cite{LY86, LiP-book, SchYau, YauH}). 
Such properties include (but are not restricted to) H\"older regularity and higher order differentiability, sharp spectral asymptotics 
and bounds, heat kernel bounds, 
Liouville-type results and many more ({\it cf.} \cite{AM, BDM, Bri, Giaq, Gr, LiX, SZ, Taheri-book-one, Taheri-book-two, Zhang}).

Whilst in proving gradient estimates one typically works with an explicit nonlinearity with a specific structure 
(of singularity, regularity, decay and growth), in this paper we keep the analysis and discussion on a fairly 
general level without confining to specific examples in order to firstly provide a unified treatment 
of the estimates and secondly to more clearly see how the structure and form of the nonlinearity influences 
the estimates and subsequent results. As such our approach and analysis largely unify, extend and at places, 
improve various existing results in the literature for specific choices of nonlinearities (see below for more).

Gradient estimates for positive solutions to \eqref{eq11} in the special case of the nonlinearity 
being a superposition of a logarithmic and a linear term with variable coefficients:   
\begin{equation} \label{eq1.4}
\Delta_f u (x) + \mathsf{p}(x) u (x) \log u(x) + \mathsf{q}(x) u(x) = 0, 
\end{equation}
along with its parabolic counterpart have been the subject of extensive studies ({\it see}, e.g., Ma \cite{Ma}, 
Ruan \cite{Ruan}, Wu \cite{ Wu10} and Yang \cite{YYY} and the references therein). The interest in such 
problems originates from its natural links with gradient Ricci solitons. Recall that a Riemannian manifold 
$(M,g)$ is said to be a gradient Ricci soliton {\it iff} there there exists a smooth function $f$ on $M$ and 
a constant $\lambda \in \mathbb{R}$ such that ({\it cf.} \cite{Cao, Chow, Lott})
\begin{equation}\label{eq1.5}
{\mathscr Ric}_f(g) = {\mathscr Ric}(g) + {\rm Hess}(f) = \lambda g.
\end{equation}
A gradient Ricci soliton can be expanding ($\lambda>0$), steady ($\lambda=0$) or shrinking ($\lambda<0$). 
The notion is a generalisation of an Einstein manifold and has a fundamental role in the analysis of singularities 
of the Ricci flow \cite{Ham, Zhang}. Taking trace from both sides of \eqref{eq1.5} and using the contracted 
Bianchi identity leads one to a simple form of \eqref{eq1.4} with constant coefficients: 
$\Delta u + 2 \lambda u \log u = (A_0 - n \lambda)u$ for suitable constant $A_0$ and $u=e^f$ 
({\it see} \cite{Ma} for details). Other types of equations closely relating to \eqref{eq1.4} including:  
\begin{equation}
\Delta_f u (x) + \mathsf{p}(x) u^a(x) |\log u|^b(x) = 0, 
\end{equation}
for real exponents $a, b$ or more generally for a nonlinear function $\gamma=\gamma(s)$ on $\mathbb{R}$: 
\begin{equation}
\Delta_f u (x) + \mathsf{p}(x) u^a(x) \gamma(\log u)(x) + \mathsf{q}(x) u^b(x) =0, 
\end{equation}
have been studied in detail in \cite{CGS, GKS, Taheri-GE-1, Taheri-GE-2, Wu10, Wu15}.

Yamabe type equations $\Delta u + \mathsf{p}(x) u^s + \mathsf{q}(x) u  =0$ are also of form  
\eqref{eq11} with a power-like nonlinearity. 
Bidaut-V\'eron and V\'eron \cite{Bid} 
studied the equation $\Delta u + u^s+ \mathsf{q} u =0$ on a compact manifold and under suitable conditions 
on the Ricci tensor, $n$ and $s, \mathsf{q}$ showed that it only admits constant 
solutions. Gidas and Spruck \cite{GidSp} considered 
\begin{equation} \label{eq1.12}
\Delta u(x) + \mathsf{p}(x) u^s(x) =0, \qquad 1 \le s < (n+2)/(n-2),
\end{equation} 
and showed that when ${\mathscr Ric}(g) \ge 0$ any non-negative solution to this equation must be zero. 
Yang \cite{Yang} showed that the same equation with constant $\mathsf{p}>0$ and $s<0$ admits no positive 
solution when ${\mathscr Ric}(g) \ge 0$. Note that the case $s=3$ in \eqref{eq1.12} is related 
to Yang-Mills equation ({\it cf.} Cafarelli, Gidas and Spruck \cite{CGS}) and the case $s<0$ is related to the 
steady states of the thin films equation ({\it cf.} Guo and Wei \cite{Guo}). For more related results see Brandolini, Rigoli and Setti \cite{BRS}, 
Li \cite{LiJ91}, Li, Tam and Yang\cite{LiPTam} and Zhang \cite{Zhang-Pac} and for a more detailed account 
on the Yamabe problem in geometry see \cite{Lee, Mast}. The natural form of Yamabe equation 
in the setting of smooth metric measure spaces is   
\begin{equation} 
\Delta_f u (x) + \mathsf{p}(x) u^s(x) + \mathsf{q}(x) u (x)  = 0. 
\end{equation}
For gradient estimates, Harnack inequalities and other counterparts of the above results we 
refer the reader to Case \cite{Case}, Wu \cite{Wu15}, Zhang and Ma \cite{ZhMa}.

A more general form of Yamabe equation is the Einstein-scalar field Lichnerowicz equation 
(see, e.g., Choquet-Bruhat \cite{CBY}, Chow \cite{Chow}, Zhang \cite{Zhang}). 
When the underlying manifold has dimension $n \ge 3$ this takes the form 
$\Delta u + \mathsf{p}(x)u^\alpha + \mathsf{q}(x) u^\beta + \mathsf{r}(x) u =0$ with 
$\alpha=(n+2)/(n-2)$ and $\beta=(3n-2)/(n-2)$ while when $n=2$ this takes the 
form $\Delta u + \mathsf{p}(x) e^{2u} + \mathsf{q}(x) e^{-2u} + \mathsf{r}(x)=0$.  
The Einstein-scalar field Lichnerowicz equation in the setting of smooth metric 
measure spaces can be further generalised and written as:  
\begin{equation}
\Delta_ f u + \mathsf{p}(x) u^\alpha + \mathsf{q}(x) u^\beta + \mathsf{r}(x) u \log u + \mathsf{h}(x) u =0, 
\end{equation}
and 
\begin{equation}
\Delta_ f u + \mathsf{p}(x) e^{2u} + \mathsf{q}(x) e^{-2u} + \mathsf{r}(x) =0. 
\end{equation}

For gradient estimates, Harnack inequalities and Liouville-type results in this and related contexts see 
Dung, Khanh and Ngo \cite{Dung}, Song and Zhao \cite{Song}, Taheri \cite{Taheri-GE-1, Taheri-GE-2}, 
Wu \cite{Wu18} and the references therein.

Let us end this introduction by briefly describing the plan of the paper. In Section \ref{sec2} we present 
the main results of the paper, namely, a local and global gradient estimate for equation \eqref{eq11}, 
followed by both local and global Harnack inequalities and a general Liouville-type result. The subsequent 
sections, namely, \ref{sec3}, \ref{Harnack-Section} and \ref{Liouville-Section} are then devoted to the 
detailed proofs respectively.

\qquad \\
{\bf Notation}
Fixing a base point $p \in M$ we denote by $d=d_p(x)$ the Riemannian distance between $x$ and 
$p$ with respect to the metric $g$ and by $r=r_p(x)$ the geodesic radial variable with origin at $p$. 
We denote by $\mathcal{B}_R(p) \subset M$ the closed geodesic ball of radius $R>0$ centred at 
$p$. When the choice of the point $p$ is clear from the context we often abbreviate and write 
$d(x)$, $r(x)$ or $\mathcal{B}_R$ respectively. We write $s_+=\max(s, 0)$ and $s_-=\min(s, 0)$ 
and so $s=s_+ + s_-$ with $s_+ \ge 0$ and $s_- \le 0$. 

For given function $\Sigma=\Sigma[x,u]$ we denote its partial derivatives by subscripts, e.g., 
$\Sigma_x$, $\Sigma_u$, {\it etc.} and we reserve the notation $\Sigma^x$ for 
the function $\Sigma[\cdot, u]$ obtained by freezing the argument $u$ and viewing it as 
a function of $x$; e.g., below we frequently use $\nabla \Sigma^x$ and $\Delta_f \Sigma^x$.   

For the sake of reader's convenience we recall that in local coordinates $(x^i)$ we have the following 
formulae for the Laplace-Beltrami operator, Riemann and Ricci curvature tensors respectively:  
\begin{equation}
\Delta = \frac{1}{\sqrt{|g|}} \frac{\partial}{\partial x_i} \left( \sqrt{|g|} g^{ij} \frac{\partial}{\partial x_j} \right),
\end{equation}
and 
 \begin{equation}
[{\mathscr Rm}(g)]^\ell_{ijk} = \frac{\partial \Gamma^\ell_{jk}}{\partial x_i} - \frac{\partial \Gamma^\ell_{ik}}{\partial x_j} 
+ \Gamma^p_{jk} \Gamma^\ell_{ip}  - \Gamma^p_{ik} \Gamma^\ell_{jp},   
\end{equation}
and 
\begin{equation}
[{\mathscr Ric}(g)]_{ij} = \frac{\partial \Gamma_{ij}^k}{\partial x_k} - \frac{\partial \Gamma_{\ell j}^\ell}{\partial x_i} 
+ \Gamma^k_{ij} \Gamma^\ell_{\ell k} - \Gamma^\ell_{ik} \Gamma^k_{\ell j}. 
\end{equation}
Note that here
\begin{equation}
\Gamma^k_{ij} = \frac{1}{2} g^{k\ell} \left( \frac{\partial g_{j \ell}}{\partial x_i} 
+ \frac{\partial g_{i \ell}}{\partial x_j} - \frac{\partial g_{ij}}{\partial x_\ell} \right),
\end{equation}
are the Christoffel symbols and $g_{ij}$, $|g|$ and $g^{ij} = (g^{-1})_{ij}$ are respectively the components, 
determinant and the components of the inverse of the metric tensor $g$.

\section{Statement of the main results} \label{sec2}

In this section we present the main results of the paper. The proofs are delegated to the subsequent sections. We 
emphasise that throughout the paper, the curvature lower bounds are expressed in the form ${\mathscr Ric}_f^m(g) \ge - (m-1)k g$ 
with $k \ge 0$ and $m$ a suitable constant $n \le m <\infty$. As the estimates here are of Li-Yau type, it is well-known 
that a lower bound on ${\mathscr Ric}_f(g)$ is not sufficient for the purpose.

\subsection{A local and a global Li-Yau type gradient estimate for \eqref{eq11}}

\begin{theorem} \label{thm1}
Let $(M, g, d\mu)$ be a complete smooth metric measure space with $d\mu=e^{-f} dv_g$ and assume that 
${\mathscr Ric}_f^m(g) \ge -(m-1) k g$ in $\mathcal{B}_{2R} = \mathcal{B}_{2R}(p)$ for suitable constants 
$m \ge n$ and $k \ge 0$.  Let $u$ be a positive solution to \eqref{eq11} in $\mathcal{B}_{2R}$. Then for 
every $\mu>1$ and $\varepsilon \in (0,1)$ and every $x \in \mathcal{B}_R$,

\begin{align}\label{1.26}
\frac{|\nabla u|^2}{\mu u^2} + \frac{\Sigma [x,u]}{u} \le&~ 
\frac{m\mu}{2R^2}\left[ (c_2 +(m-1) c_1(1+R \sqrt{k})+2c_1^2) + \frac{mc_1^2 \mu^2}{4(\mu-1)} \right] \nonumber\\
&+ \frac{\sqrt{m}}{2} \bigg[ 
\frac{m \mu^2 \mathsf{A}_\Sigma^2}{(1-\varepsilon)(\mu-1)^2} 
+ \left[\frac{27m \mu^2 \mathsf{B}_\Sigma^4}{4\varepsilon(\mu-1)^2}\right]^{1/3} 
+ 2\mu \mathsf{C}_\Sigma \bigg]^{1/2} \nonumber \\
&+ \frac{m \mu}{2} \sup_{\mathcal{B}_{2R}} \left\{ \frac{(u \Sigma_u[x,u]-\Sigma[x,u] )_+}{u} \right\},
\end{align}
where the quantities $\mathsf{A}_\Sigma$, $\mathsf{B}_\Sigma$ and $\mathsf{C}_\Sigma$ are given by 
\begin{align} \label{eq2.2}
\mathsf{A}_\Sigma = \sup_{\mathcal{B}_{2R}} \left\{ \frac{2(m-1)k u 
+ (-\Sigma[x,u] + u\Sigma_u[x,u] - \mu u^2 \Sigma_{uu}[x,u])_+}{2u} \right\}, 
\end{align}
\begin{align} \label{eq2.3}
\mathsf{B}_\Sigma = \sup_{\mathcal{B}_{2R}} \left\{ \frac{|\Sigma_x[x,u] -\mu u \Sigma_{xu}[x,u]|}{u} \right\}, \qquad 
\mathsf{C}_\Sigma = \sup _{\mathcal{B}_{2R}} \left\{ \frac {(-\Delta_f \Sigma^x [x,u])_+}{u} \right\}.
\end{align}
\end{theorem}

The local estimate above has a global counterpart subject to the prescribed bounds in the theorem being global. 
The proof follows by passing to the limit $R \to \infty$ in \eqref{1.26} and taking into account the vanishing of certain 
terms as a result of the bounds being global and the relevant constants being independent of $R$. The precise 
formulation of this is given in the following theorem.

\begin{theorem} \label{thm1-global}
Let $(M, g, d\mu)$ be a complete smooth metric measure space with $d\mu=e^{-f} dv_g$ and 
${\mathscr Ric}_f^m(g) \ge -(m-1) k g$ on $M$ where $m \ge n$ and $k \ge 0$. 
Assume $u$ is a positive solution to \eqref{eq11} on $M$. Then for every $\mu>1$ 
and $\varepsilon \in (0,1)$ and every $x \in M$,

\begin{align}\label{1.26-global}
\frac{|\nabla u|^2}{\mu u^2} + \frac{\Sigma[x,u]}{u} \le&~ 
\frac{\sqrt{m}}{2} \bigg[ \frac{m \mu^2\mathsf{A}_\Sigma^2}{(1-\varepsilon)(\mu-1)^2} 
+ \left[\frac{27 m \mu^2\mathsf{B}_\Sigma^4}{4 \varepsilon(\mu-1)^2}\right]^{1/3}
+ 2 \mu \mathsf{C}_\Sigma \bigg]^{1/2} \nonumber\\ 
& + \frac{m \mu}{2} \sup_{M}\left\{ \frac{(u\Sigma_u[x,u]-\Sigma[x,u])_+}{u} \right \}.
\end{align}
Here $\mathsf{A}_\Sigma$, $\mathsf{B}_\Sigma$ and $\mathsf{C}_\Sigma$ are as in \eqref{eq2.2}-\eqref{eq2.3} 
in Theorem $\ref{thm1}$ except that now the supremums are taken over all of $M$.
\end{theorem}

\subsection{A local and a global elliptic Harnack inequality for \eqref{eq11}}

\begin{theorem} \label{cor Harnack}
Under the assumptions of Theorem $\ref{thm1}$ and the Bakry-\'Emery curvature 
bound ${\mathscr Ric}_f^m(g) \ge - (m-1) k g$ on $\mathcal{B}_{2R}$, for any 
positive solution $u$ to \eqref{eq11}, and any $x_1, x_2 \in \mathcal{B}_R$ we have 
\begin{equation} \label{eq2.5}
u(x_2) \le e^{2R \sqrt \mathbb{H} } u(x_1). 
\end{equation}
The positive constant $\mathbb{H}$ can be explicitly expressed in terms of the local bounds as 
\begin{align}
\mathbb{H} = &~m\mu^2[ (c_2 +(m-1) c_1(1+R \sqrt{k})+2c_1^2) + mc_1^2 \mu^2/(4(\mu-1))] /(2R^2)\nonumber\\
&+\sqrt {m} \mu/2\{ m \mu^2\mathsf{A}_\Sigma^2/[(1-\varepsilon)(\mu-1)^2] 
+ [27m \mu^2\mathsf{B}_\Sigma^4/ (4\varepsilon(\mu-1)^2)]^{1/3} 
+ 2\mu \mathsf{C}_\Sigma \} ^{1/2} \nonumber \\
&+ m \mu^2 \sup_{\mathcal{B}_{2R}} \{{(u\Sigma_u[x,u]-\Sigma[x,u])_+/(2u)}\} \nonumber \\
&- \mu \inf_{\mathcal{B}_{R}} \{{(\Sigma[x,u])_-/u}\}, 
\end{align}
where $\mathsf{A}_\Sigma$, $\mathsf{B}_\Sigma$ and $\mathsf{C}_\Sigma$ are as in \eqref{eq2.2}-\eqref{eq2.3} 
in Theorem $\ref{thm1}$. In particular from \eqref{eq2.5} we have 
\begin{equation}
\sup_{\mathcal{B}_R} u \le e^{2R \sqrt \mathbb{H}}  \inf_{\mathcal{B}_R} u.
\end{equation}
\end{theorem}

For the global version we can use a similar argument utilising the global bounds in Theorem \ref{thm1-global} 
and have the counterpart of \eqref{eq2.5} with $d(x_1, x_2)$ replacing $2R$ and $\mathbb{H}>0$ now being 
its global version from \eqref{1.26-global}. The precise formulation is given below.

\begin{theorem} \label{cor Harnack global}
Under the assumptions of Theorem $\ref{thm1-global}$ and the Bakry-\'Emery curvature 
bound ${\mathscr Ric}_f^m(g) \ge - (m-1) k g$ on $M$, for any 
positive solution $u$ to \eqref{eq11}, and any $x_1, x_2 \in M$ we have 
\begin{equation} \label{eq2.5-global}
u(x_2) \le e^{d(x_1, x_2) \sqrt \mathbb{H}} u(x_1). 
\end{equation}
The positive constant $\mathbb{H}$ can be explicitly expressed in terms of the global bounds as 
\begin{align}
\mathbb{H} = &~\sqrt {m} \mu/2\{ m \mu^2\mathsf{A}_\Sigma^2/[(1-\varepsilon)(\mu-1)^2] 
+ [27m \mu^2\mathsf{B}_\Sigma^4/ (4\varepsilon(\mu-1)^2)]^{1/3} 
+ 2\mu \mathsf{C}_\Sigma \} ^{1/2} \nonumber \\
&+ m \mu^2 \sup_{M} \{{(u\Sigma_u[x,u]-\Sigma[x,u])_+/(2u)}\} 
- \mu \inf_{M} \{{(\Sigma[x,u])_-/u}\}, 
\end{align}
where $\mathsf{A}_\Sigma$, $\mathsf{B}_\Sigma$ and $\mathsf{C}_\Sigma$ are as in Theorem $\ref{thm1-global}$.  
\end{theorem}

\subsection{A Liouville-type theorem and some applications}

\begin{theorem} \label{thm Liouville}
Let $(M,g, d\mu)$ be a smooth metric measure space with $d\mu=e^{-f} dv_g$ satisfying ${\mathscr Ric}^m_f(g) \ge 0$ in $M$. 
Let $u$ be a positive solution to $\Delta_f u + \Sigma[u] =0$. Assume that $\Sigma[u] \ge 0$, $u\Sigma_u [u] -\Sigma[u] \le 0$ 
and $\mu u^2 \Sigma_{uu} [u]  - u \Sigma_u [u] + \Sigma [u] \ge 0$ for some $\mu >1$. Then $u$ must be a constant. 
In particular $\Sigma[u] =0$.
\end{theorem}

The proof of this theorem follows from the gradient estimates established above and is presented 
in Section \ref{Liouville-Section}. We end this section by giving two stark applications of the above 
theorem. To this end consider first a superposition of power-like nonlinearities with real coefficients 
$\mathsf{p}_j$ and exponents $a_j$ for $1 \le j \le N$ in the form 
\begin{equation} \label{polies PR}
\Sigma[u] = \sum_{j=1}^N \mathsf{p}_j u^{a_j}.
\end{equation}
A direct calculation gives  
\begin{align}
u \Sigma_u [u] - \Sigma [u] &= \sum_{j=1}^N \mathsf{p}_j (a_j -1) u^{a_j}, \\
\mu u^2 \Sigma_{uu} [u]  - u \Sigma_u [u] + \Sigma [u] 
&= \sum_{j=1}^N [\mu \mathsf{p}_j a_j (a_j -1) - \mathsf{p}_j a_j + \mathsf{p}_j] u^{a_j} \nonumber \\
&= \sum_{j=1}^N [\mathsf{p}_j (a_j -1) (\mu a_j -1)] u^{a_j}. 
\end{align}

Evidently for the range $\mathsf{p}_j \ge 0$ we have $\Sigma(u) \ge 0$ whilst subject to $a_j \le 1$ we have 
$u\Sigma_u [u] -\Sigma [u] \le 0$ and $\mu u^2 \Sigma_{uu} [u]  - u \Sigma_u [u] + \Sigma [u] \ge 0$ (by 
choosing $\mu>1$ suitably). Theorem \ref{thm Liouville} now leads to the following conclusion extending earlier 
results on Yamabe type problems to more general nonlinearities ({\it cf.} \cite{Dung, GidSp, Yang, Wu18}). 
Further applications and results in this direction will be discussed in a forthcoming paper 
({\it see} also \cite{Taheri-GE-1, Taheri-GE-2}).

\begin{theorem} \label{LiouvilleThmEx}
Let $(M, g, d\mu)$ be a complete smooth metric measure space with $d\mu=e^{-f} dv_g$ and 
${\mathscr Ric}^m_f(g) \ge 0$. Let $u$ be a positive smooth solution to the equation   
\begin{equation} \label{elliptic PDE 2}
\Delta_f u + \sum_{j=1}^N \mathsf{p}_j u^{a_j} = 0.
\end{equation}
If $\mathsf{p}_j \ge 0$ and $a_j \le 1$ for $1 \le j \le N$ then $u$ must be a constant. 
\end{theorem}

\begin{remark} Note that a constant solution to \eqref{eq11} must be a zero of $\Sigma$. 
Thus if $\Sigma$ has no positive zeros  then the above Liouville theorem becomes a 
non-existence result. In Theorem \ref{LiouvilleThmEx} and for \eqref{polies PR} this 
happens when $\mathsf{p}_j>0$ for at least one $1 \le j \le N$.
\end{remark}

As another application, again relating to the discussions in Section \ref{sec1}, consider a superposition of 
a logarithmic and a power-like nonlinearity with real coefficients $\mathsf{p}, \mathsf{q}$, exponent $s$ 
and $\gamma \in \mathscr{C}^2(\mathbb R)$ in the form 
\begin{equation}
\Sigma[u] = \mathsf{p} u \gamma(\log u) + \mathsf{q} u^s. 
\end{equation}
A straightforward calculation then gives $u \Sigma_u [u] - \Sigma [u] = \mathsf{p} u \gamma'(\log u) + \mathsf{q} (s -1) u^{s}$ and 
$\mu u^2 \Sigma_{uu} [u]  - u \Sigma_u [u] + \Sigma [u] = {\mathsf p}[(\mu-1) u \gamma' + \mu u \gamma''] + [\mathsf{q} (s -1) (\mu s -1)] u^{s}$. 
The following theorem now directly results from Theorem \ref{thm Liouville}.

\begin{theorem} \label{LiouvilleThmEx-log}
Let $(M, g, d\mu)$ be a complete smooth metric measure space with $d\mu=e^{-f} dv_g$ and 
${\mathscr Ric}^m_f(g) \ge 0$. Let $u$ be a positive smooth solution to the equation   
\begin{equation} \label{elliptic PDE 2}
\Delta_f u + \mathsf{p} u \gamma(\log u) + \mathsf{q} u^s = 0.
\end{equation}
Assume that $\mathsf{p}, \mathsf{q} \ge 0$, $s \le 1$ and that along the solution $u$ we have 
$\gamma \ge 0$, $\gamma' \le 0$ and $\mu\gamma''+(\mu-1)\gamma'\ge0$ for some $\mu >1$ 
$($with $1<\mu<1/s$ if $0<s<1$$)$. Then $u$ must be a constant.
\end{theorem}

\section{Proof of the Li-Yau type gradient estimate in Theorem \ref{thm1}} \label{sec3}

This section is devoted to the proof of the main estimate for the positive solutions of \eqref{eq11} in its 
local form. In its global form, the estimate, as seen, follows by passing to the limit $R \to \infty$. As this 
requires a number of technical lemmas and tools, we pause briefly, to present these results and tools 
in the next subsection, before moving on to the proof of the local estimate in Theorem \ref{thm1} in 
the following subsection.

\subsection{Some technical lemmas and identities}

\begin{lemma}
Let $u$ be a positive solution to the equation \eqref{eq11} and let $h = \log u$. Then $h$ satisfies the equation 
\begin{align}\label{3.1}
\Delta_f h + |\nabla h|^2 + e^{-h} \Sigma [x,e^h] =0.
\end{align}
\end{lemma}

\begin{proof}
An easy calculation gives $\nabla h = (\nabla u)/u$ and $\Delta h=(\Delta u)/u - |\nabla u|^2/u^2$. 
Hence $\Delta_f h = (\Delta_f u)/u - |\nabla u|^2/u^2=-\Sigma [x,u]/u - |\nabla u|^2/u^2$ giving the 
desired conclusion. 
\end{proof}

\begin {lemma}\label{H-fLap-identity}
Let $u$ be a positive solution to \eqref{eq11}, $h=\log u$ and let $H$ be defined by 
\begin{align} \label{3.2}
H= |\nabla h|^2 + \mu e^{-h} \Sigma [x,e^h], 
\end{align}
where $\mu \ge 1$ is an arbitrary constant. Then $H$ satisfies the equation 
\begin{align} \label{eq3.2}
\Delta _f H =&~ 2|\nabla^2 h|^2 +2 \frac{\langle \nabla f , \nabla h \rangle^2}{m-n}- 2 \langle \nabla h , \nabla H \rangle
+ 2 {\mathscr Ric}_f^m (\nabla h, \nabla h) \nonumber\\
& - 2(\mu -1) e^{-h}\Sigma |\nabla h|^2 + 2(\mu -1) e^{-h} \langle \nabla h,\nabla \Sigma \rangle + \mu \Delta_f (e^{-h}\Sigma), 
\end{align}
where we have abbreviated the arguments of $\Sigma$ and its derivatives. 
\end{lemma}

\begin{proof}
Referring to the formulation of $H$ is \eqref{3.2} an application of $\Delta_f$ to both sides of the equation gives 
\begin{align}\label{3.5}
\Delta_f H = \Delta_f |\nabla h|^2 + \mu \Delta_f (e^{-h} \Sigma [x,e^h]).
\end{align}
Furthermore, referring to \eqref{3.1} and again to \eqref{3.2} it is evident that we have the relation:  
\begin{align}\label{3.4}
\Delta_f h = - (|\nabla h|^2 + e^{-h} \Sigma [x, e^h]) = -(H - (\mu -1) e^{-h} \Sigma [x,e^h]).
\end{align}

Now as for the first term on the right-hand side of \eqref{3.5}, by the generalised Bochner-Weitzenb\"ock formula (as applied to $h$), we have
\begin {align} \label{3.6}
\Delta_f |\nabla h|^2 = 2|\nabla^2 h|^2 + 2\langle \nabla h , \nabla \Delta_f h \rangle 
+ 2{\mathscr Ric}_f^m (\nabla h, \nabla h) 
+2 \frac{\langle \nabla f , \nabla h\rangle^2}{m-n}.  
\end{align}  
Hence, by substituting \eqref{3.6} into \eqref{3.5} and making note of \eqref{3.4}, we have after a basic differentiation, 
\begin{align}\label{3.7}
\Delta_f H = &~ 2|\nabla^2 h|^2 + 2\langle \nabla h , \nabla (-H +(\mu -1) e^{-h}\Sigma [x,e^h]) \rangle \nonumber\\
&+ 2{\mathscr Ric}_f^m (\nabla h, \nabla h) +2 \frac{\langle \nabla f , \nabla h\rangle^2}{m-n} 
+ \mu \Delta_f (e^{-f}\Sigma [x,e^h]) \nonumber\\
=&~ 2|\nabla^2 h|^2 - 2\langle \nabla h , \nabla H \rangle - 2 (\mu -1) e^{-h} \Sigma [x,e^h] |\nabla h|^2 \nonumber\\
& + 2 (\mu -1) e^{-h}\langle \nabla h , \nabla \Sigma [x,e^h] \rangle + 2{\mathscr Ric}_f^m (\nabla h, \nabla h)\nonumber\\ 
&+2 \frac{\langle \nabla f , \nabla h\rangle^2}{m-n} + \mu \Delta_f (e^{-h}\Sigma [x,e^h]), 
\end{align}
which upon a rearrangement of terms gives the desired identity. 
\end{proof}

\begin{lemma} \label{lem3.3}
Let $u$ be a positive solution to \eqref{eq11}, $h=\log u$ and let $H$ be as defined by \eqref{3.2}. Then, if 
${\mathscr Ric}_f^m(g) \geq -(m-1) kg$, we have 
\begin{align}
\Delta_f H \ge&~ 2 \frac{(\Delta_f h)^2}{m} - 2 \langle \nabla h , \nabla H \rangle + [e^{-h}\Sigma - \Sigma_u] H\nonumber\\
&+\left[ e^{-h} \Sigma - \Sigma_u + \mu e^h \Sigma_{uu} -2(m-1)k \right] |\nabla h|^2 \nonumber\\
& -2 \langle \nabla h , [e^{-h} \Sigma_x - \mu \Sigma_{xu}] \rangle + \mu e^{-h} \Delta_f \Sigma^{x}.
\end{align}
\end{lemma}

\begin{proof} 
By virtue of the bound ${\mathscr Ric}_f^m(g) \ge -(m-1) k g$ it follows upon recalling the identity in Lemma \ref{H-fLap-identity} that 
\begin{align} \label {3.16}
\Delta_f H \ge&~ 2 \frac{(\Delta_f h)^2}{m}- 2\langle \nabla h , \nabla H \rangle - 2 (\mu -1) e^{-h} \Sigma |\nabla h|^2 \nonumber\\
&+ 2 (\mu -1) e^{-h}\langle \nabla h , \nabla \Sigma \rangle - 2 (m-1) k |\nabla h|^2\nonumber\\ 
& + \mu \Sigma \Delta_f e^{-h}+2\mu \langle \nabla e^{-h} , \nabla \Sigma \rangle + \mu e^{-h} \Delta_f \Sigma.
\end{align}

Note that in concluding the above inequality, specifically, the first term on the right-hand side, we have made use of the basic inequalities  
\begin {align} \label {3.15} 
|\nabla^2 h|^2+ \frac{\langle \nabla f , \nabla h \rangle^2}{m-n} 
\ge \frac{(\Delta h)^2}{n} + \frac{\langle \nabla f, \nabla h \rangle ^2}{m-n}\ge \frac{(\Delta_f h)^2}{m}. 
\end{align}

Let us now proceed by attending to some useful and straightforward calculations relating to the nonlinear term $\Sigma=\Sigma(x, e^h)$. 
Evidently 
\begin {align} \label{3.8}
\nabla \Sigma [x,e^h] = \Sigma_x [x,e^h] +e^h \Sigma_u [x,e^h] \nabla h,
\end{align}  
and thus moving on to the Laplacian we can write 
\begin {align} 
\Delta \Sigma [x,e^h] &= {\rm div} ( \Sigma_x [x,e^h] + e^h \Sigma_u [x,e^h] \nabla h ). 
\end{align}
It is convenient to do the calculations in local coordinates and so we proceed by writing 
\begin{align}
\Delta \Sigma (x,e^h) 
=& \sum_{i=1}^n \frac{\partial}{\partial x_i} ( \Sigma_{x_i} [x,e^h] + e^h \Sigma_u [x,e^h] h_i ) \nonumber\\
=& \sum_{i=1}^n \bigg\{ \Sigma_{x_i x_i} [x,e^h] + \Sigma _{x_i u} [x,e^h] (e^h)_i + e^h h_i \Sigma_u [x,e^h] h_i \nonumber \\
& + e^h (\Sigma_{x_i u} [x, e^h] + \Sigma_{uu} [x,e^h] (e^h)_i) h_i 
+ e^h \Sigma_u [x, e^h] h_{ii} \bigg\}.  
\end{align}
Abbreviating the arguments $[x, e^h]$ of $\Sigma$ for convenience and rewriting the above we have
\begin{align}
\Delta \Sigma
&= \Delta \Sigma^{x}+e^h \langle \Sigma_{xu} ,\nabla h \rangle +e^h |\nabla h|^2 \Sigma_u \nonumber \\
& +e^h \langle \Sigma_{xu},\nabla h \rangle+e^{2h} |\nabla h|^2 \Sigma_{uu} +e^h \Sigma_u \Delta h \nonumber\\
& = \Delta \Sigma^{x}+2e^h \langle \Sigma_{xu} ,\nabla h \rangle + e^h |\nabla h|^2 ( \Sigma_u +e^h \Sigma_{uu})+e^h \Sigma_u \Delta h.
\end{align}   
As a result the above upon substitution give
\begin {align} \label{3.12}
\Delta_f \Sigma =&~ \Delta \Sigma -\langle \nabla f , \nabla \Sigma \rangle = \Delta \Sigma 
-\langle \nabla f , ( \Sigma_x+e^h \Sigma_u\nabla h) \rangle \nonumber\\
=&~\Delta \Sigma -\langle \nabla f, \Sigma_x \rangle -e^h \Sigma_u \langle \nabla f , \nabla h \rangle \nonumber\\
=&~\Delta \Sigma^{x} -\langle \nabla f, \nabla \Sigma^x \rangle +2e^h \langle \Sigma_{xu} ,\nabla h \rangle \nonumber \\
&+ e^h |\nabla h|^2 ( \Sigma_u +e^h \Sigma_{uu})+e^h \Sigma_u \Delta h - e^h \Sigma_u \langle \nabla f, \nabla h \rangle \nonumber\\
=&~\Delta_f \Sigma^{x} +2e^h \langle \Sigma_{xu} ,\nabla h \rangle 
+ e^h |\nabla h|^2 ( \Sigma_u +e^h \Sigma_{uu})+e^h \Sigma_u \Delta_f h.
\end{align}
Moreover for the sake of future reference we also note that 
\begin {align} \label{3.13}
\Delta_f e^{-h} &= \Delta e^{-h} - \langle \nabla f , \nabla e^{-h} \rangle \nonumber \\
&= {\rm div}(-e^{-h}\nabla h)+e^{-h} \langle \nabla f, \nabla h \rangle \nonumber\\
&= -e^{-h} \Delta h + e^{-h} |\nabla h|^2 +e^{-h} \langle \nabla f , \nabla h \rangle \nonumber \\
&= -e^{-h} (\Delta_f h - |\nabla h|^2).
\end{align}

Now returning to the inequality \eqref{3.16} and upon substituting from \eqref{3.8}, \eqref{3.12} and \eqref {3.13} we obtain  
\begin{align}
\Delta_f H 
\ge&~2(\Delta_f h)^2/m - 2\langle \nabla h , \nabla H \rangle - 2 (\mu -1) e^{-h} \Sigma |\nabla h|^2 \nonumber\\
& + 2 (\mu -1) e^{-h}\langle \nabla h , \Sigma_x +e^h \Sigma_u \nabla h \rangle - 2 (m-1) k |\nabla h|^2\nonumber\\ 
& - \mu e^{-h} \Sigma ( \Delta_f h - |\nabla h|^2 )  -2 \mu e^{-h} \langle \nabla h, \Sigma_x +e^h \Sigma_u \nabla h \rangle \nonumber\\
& +\mu e^{-h} [\Delta_f \Sigma^{x} +2e^h \langle \Sigma_{xu} ,\nabla h \rangle 
+ e^h ( \Sigma_u +e^h \Sigma_{uu}) |\nabla h|^2 + e^h \Sigma_u \Delta_f h], 
\end{align}
or upon rearranging 
\begin{align} \label{3.17}
\Delta_f H 
\ge&~2 (\Delta_f h)^2/m- 2\langle \nabla h , \nabla H \rangle - 2 (\mu -1) e^{-h} \Sigma |\nabla h|^2 \nonumber\\
&+ 2 (\mu -1) e^{-h}\langle \nabla h , \Sigma_x \rangle + 2 (\mu -1) \Sigma_u |\nabla h|^2 \nonumber \\
&- 2 (m-1) k |\nabla h|^2 - \mu (e^{-h} \Sigma -\Sigma_u) \Delta_f h + \mu e^{-h} \Sigma |\nabla h|^2 \nonumber \\
& - 2 \mu e^{-h} \langle \nabla h,\Sigma_x \rangle 
- 2 \mu \Sigma_u |\nabla h|^2  +\mu e^{-h} \Delta_f \Sigma^{x} \nonumber \\
& +2\mu \langle \nabla h ,\Sigma_{xu} \rangle 
+ \mu ( \Sigma_u +e^h \Sigma_{uu})  |\nabla h|^2.
\end{align}
Next by recalling \eqref{3.2} and \eqref{3.4} we can write
\begin{align}\label{3.18}
\mu \Delta_f h = - \mu (H - (\mu -1) e^{-h}\Sigma [x,e^h]) =- [H + (\mu -1) |\nabla h|^2], 
\end{align}
and therefore by substituting the latter back in \eqref{3.17} and rearranging terms it follows that 
\begin{align}\label{3.19}
\Delta_f H \ge&~2 (\Delta_f h)^2/m- 2\langle \nabla h , \nabla H \rangle \nonumber \\
& + (e^{-h}\Sigma -\Sigma_u) [H+(\mu -1)|\nabla h|^2] \nonumber\\
& + [-2(\mu -1) e^{-h}\Sigma -2 \Sigma_u -2(m-1) k] |\nabla h|^2 \nonumber \\
& + (\mu e^{-h} \Sigma + \mu \Sigma_u + \mu e^h \Sigma_{uu}) |\nabla h|^2 \nonumber\\
& -2e^{-h} \langle \nabla h , \Sigma_x \rangle + \mu e^{-h} \Delta_f \Sigma^{x}+2\mu \langle \nabla h , \Sigma_{xu} \rangle.
\end{align}
Finally taking into account the necessary cancellations and by a further rearrangement of terms we obtain  
\begin{align}
\Delta_f H \ge&~ 2 (\Delta_f h)^2/m - 2\langle \nabla h , \nabla H \rangle + (e^{-h}\Sigma -\Sigma_u) H  \nonumber\\
&+ [e^{-h}\Sigma -\Sigma_u + \mu e^h \Sigma_{uu} -2(m-1)k] |\nabla h|^2 \nonumber\\
&-2\langle \nabla h , e^{-h} \Sigma_x - \mu \Sigma_{xu} \rangle + \mu e^{-h}\Delta_f \Sigma^{x}, 
\end{align}
which is the desired conclusion. 
\end{proof}

The following lemma will also be used in the course of the proof of the local estimate in the next subsection.

\begin{lemma} \label{LiYau-basiclemma}
Suppose $a,b,z \in \mathbb{R}$, $c, y>0$ and $\mu>1$ are arbitrary constants such that $y-\mu z >0$. Then for any 
$\varepsilon \in (0,1)$ we have
\begin{align} \label{eq6.9} 
(y-z)^2& - a \sqrt y (y-\mu z) - b y - c \sqrt y \nonumber\\
&\ge (y-\mu z)^2/\mu^2-a^2\mu^2 (y-\mu z)/[8 (\mu-1)] \nonumber \\
&-(3/4) c^{4/3} [\mu^2/(4 \varepsilon (\mu -1)^2)]^{1/3} 
- (\mu^2 b^2)/[4(1-\varepsilon)(\mu-1)^2].  
\end{align}
\end{lemma}

\begin{proof}
Starting from the expression on the left-hand side in \eqref{eq6.9} we can write for any $\delta, \varepsilon$ by basic considerations
\begin {align} \label{eq6.10}
(y-z)^2 &- a \sqrt y (y-\mu z) - b y - c \sqrt y\nonumber\\
=&~(1-\varepsilon-\delta)y^2-(2-\varepsilon \mu)yz+z^2+(\varepsilon y - a \sqrt y)(y-\mu z) 
+ \delta y^2 - by - c \sqrt y \nonumber\\
=&~(1/\mu - \varepsilon/2)(y-\mu z)^2+(1-\varepsilon-\delta-1/\mu+\varepsilon/2)y^2
+(1- \mu + \varepsilon \mu^2/2)z^2\nonumber\\
&+ (\varepsilon y - a \sqrt y)(y-\mu z) + \delta y^2 - by - c \sqrt y. 
\end{align}   
In particular setting $\delta = (1/\mu-1)^2$ and $\varepsilon = 2-2/\mu-2(1/\mu-1)^2 = 2(\mu-1)/\mu^2$ 
gives $1-\varepsilon-\delta-1/\mu+\varepsilon/2=0$ and $1-\mu + \varepsilon\mu^2/2=0$ and so by making note 
of the inequality $\varepsilon y - a \sqrt y \ge - a^2/(4 \varepsilon)$ with $\varepsilon=2(\mu-1)/\mu^2>0$ we can 
deduce from 
\eqref{eq6.10} that
\begin {align}\label{8.36}   
(y-z)^2& - a \sqrt y (y-\mu z) - b y - c \sqrt y \\
&\ge (y-\mu z)^2/\mu^2-a^2\lambda^2 (y-\mu z)/[8 (\mu-1)] 
+ (\mu -1)^2 y^2/\mu^2 - by- c \sqrt y. \nonumber 
\end{align} 
Next, considering the last three terms only we can write, for any $\varepsilon \in (0,1)$,
\begin {align} 
(\mu -1)^2 &y^2/\mu^2 - by - c \sqrt y \nonumber\\
&\ge (\mu -1)^2y^2/\mu^2 - (1-\varepsilon)(\mu-1)^2 y^2/\mu^2 - 
(\mu^2 b^2)/[4(1-\varepsilon)(\mu-1)^2] - c \sqrt y\nonumber\\
&\ge \varepsilon (\mu -1)^2y^2/\mu^2 - 
(\mu^2 b^2)/[4(1-\varepsilon)(\mu-1)^2] - c \sqrt y\nonumber\\
&\ge -(3/4) c^{4/3} [\mu^2/(4 \varepsilon (\mu -1)^2)]^{1/3} 
- (\mu^2 b^2)/[4(1-\varepsilon)(\mu-1)^2] 
\end{align} 
where above we have made use of $(1-\varepsilon)(\mu-1)^2 y^2/\mu^2 - by \ge -(\mu^2 b^2)/[4(1-\varepsilon)(\mu-1)^2]$ 
and $\varepsilon (\mu -1)^2y^2/\mu^2 - c \sqrt y \ge  -(3/4) c^{4/3} [\mu^2/(4 \varepsilon (\mu -1)^2)]^{1/3}$  
to deduce the first and last inequalities respectively. Substituting back in \eqref{8.36} gives the desired inequality.  
\end{proof}

\subsection{Proof of Theorem \ref{thm1}}

This is based on the estimate established in Lemme \ref{lem3.3} and a localisation argument. In order to carry out the localisation 
and the relevant estimates we proceed by first constructing a suitable cut-off function. Towards this end we begin by introducing a 
function $\bar{\psi} = \bar{\psi}(t)$ satisfying the following conditions: 
\begin{enumerate}[label=$(\roman*)$]
\item $\bar\psi$ is of class $\mathscr{C}^2 [0, \infty)$.
\item $0 \le \bar\psi(t) \le 1$ for $0 \le t < \infty$ with 
\begin{align}
\bar{\psi}(t) 
= \left\{ 
\begin{array}{rcl}
1& & t \le 1, \\ 
0 & & t \ge 2.  
\end{array}
\right.
\end{align}
\item $\bar\psi$ is non-increasing (that is, $\bar\psi' \le 0$), and additionally, for suitable constants $c_1, c_2>0$, its first and second 
order derivatives satisfy the bounds 
\begin{align} \label{9.30}
- c_1 \le \frac{\bar{\psi}^{'}}{\sqrt{ \bar{\psi}}} \leq 0 \qquad and \qquad  \bar{\psi}^{''} \ge -c_2. 
\end{align}
\end{enumerate}
Next pick and fix a reference point $p$ in $M$ and with $r=r_p(x)$ set
\begin{align}\label{7.19}
\psi(x)=\bar \psi \left( \frac{r(x)}{R} \right).
\end{align}

It is evident that $\psi \equiv 1$ for when $0 \le r(x) \le R$ and $\psi \equiv 0$ for 
when $r(x) \ge 2R$. Let us now consider the localised function $\psi H$ supported on $\mathcal{B}_{2R}$. Let us also assume that 
$x_1$ in $\mathcal{B}_{2R}$ is the point where $\psi H$ attains its maximum. As for $\psi H \le 0$ the estimate 
is trivially true we can assume that $[\psi H](x_1)>0$. Furthermore by an argument of Calabi \cite{Calabi} we can assume 
that $x_1$ is not in the cut locus of $p$. It then follows that at this point:
\begin{align} \label{9.32}
\nabla (\psi H)=0 \qquad  and \qquad \Delta(\psi H) \leq 0,  \qquad \Delta_f (\psi H) \le 0.
\end{align}

Now starting from the basic identity
\begin{align}
\Delta_f (\psi H) = H \Delta_f \psi + 2 \langle \nabla \psi, \nabla H \rangle + \psi \Delta_f H
\end{align}
and making note of the relations \eqref{9.32} at the maximum point $x_1$, we can write 
\begin{align} \label{7.26}
0 &\ge H \Delta_f \psi + 2 \langle \nabla \psi, \nabla H \rangle + \psi \Delta_f H \nonumber \\
& \ge H \Delta_f \psi + \frac{2}{\psi} \langle \nabla \psi , \nabla (\psi H) \rangle -2 \frac{|\nabla \psi |^2}{\psi}H + \psi \Delta_f H \nonumber \\
& \ge H \Delta_f \psi -2 \frac{|\nabla \psi |^2}{\psi}H + \psi \Delta_f H.
\end{align}

We proceed now by obtaining suitable lower bounds for each of the three individual terms on the right-hand side of this inequality. 
As for the first term, referring to \eqref{7.19}, a straightforward calculation gives $\nabla \psi = (\bar\psi'/R) \nabla r$ and 
$\Delta \psi = \bar\psi'' |\nabla r|^2/R^2+\bar\psi' \Delta r/R$. Subsequently from these we have
\begin{equation} \label{9.35}
\Delta_f \psi = \Delta \psi - \langle \nabla f, \nabla \psi \rangle = \frac{\bar\psi''}{R^2} |\nabla r|^2+\frac{\bar\psi'}{R} \Delta_f r.
\end{equation}

For the last term on the right here the Wei-Wiley Laplacian comparison theorem ({\it cf}. \cite{[WeW09]}) together with 
${\mathscr Ric}_f^m (g) \ge -(m-1) k$ gives $\Delta_f r \le (m-1) \sqrt {k}\coth (\sqrt{k} r)$. Hence substituting back in 
\eqref{9.35} and noting $\bar\psi' \le 0$ we have:
\begin{align} \label{eq3.32}
\Delta_f \psi \ge \frac{1}{R^2} \bar\psi''+\frac {(m-1)}{R} \bar\psi' \sqrt{k} \coth (\sqrt{k} r).
\end{align}

Moreover upon noting $\coth (\sqrt{k} r) \le \coth (\sqrt{k} R)$ and $\sqrt{k} \coth (\sqrt{k} R) \le (1 + \sqrt{k} R)/R$, 
subject to $R \le r \le 2R$ [here we are using the monotonicity of $\coth x$ and the bound $x \coth x \le 1+x$ for $x>0$] 
we deduce that 
\begin{align}
(m-1) \bar\psi' \sqrt{k} \coth (\sqrt{k}r)\ge (m-1) (1+\sqrt{k}R) \bar\psi'/R.
\end{align}
Hence substituting this back in \eqref{eq3.32} and making note of the assumptions on $\bar\psi$, specifically, $0 \le \bar\psi \le 1$ 
and the lower bounds on $\bar\psi'$, $\bar\psi''$ in \eqref{9.30}, it follows that 
\begin{align} \label{9.38}
\Delta_f \psi &\ge \frac {1}{R^2}\bar{\psi}^{''}+\frac{(m-1)}{R} \left( \frac{1}{R}+\sqrt{k} \right) \bar{\psi}^{'}\nonumber\\
&\ge -\frac{c_2}{R^2}-\frac{(m-1)}{R} c_1 \left(\frac{1}{R}+\sqrt{k}\right)\nonumber\\
& =  -\frac{1}{R^2}[c_2 +(m-1) c_1(1+R \sqrt{k})], 
\end{align}
which can be readily utilised to bound the first term on the right-hand side in \eqref{7.26}. Referring next to the middle term in the 
same inequality, by using the imposed bounds on $\bar\psi'$ [the first inequality in \eqref{9.30}], we have
\begin{align}
\frac{|\nabla \psi |^2}{\psi} = \frac{\bar\psi'^2}{\bar\psi} \frac{|\nabla \varrho|^2}{R^2} 
= \left( \frac{\bar\psi'}{\sqrt {\bar\psi}} \right)^2 \frac{|\nabla \varrho|^2}{R^2} \le  \frac{c_1 ^2}{R^2}.  
\end{align}
Since for the third term on the right-hand side in \eqref{7.26} we already have the conclusion of Lemma \ref{lem3.3} at our disposal, 
by substituting the above fragments back, we obtain at the maximum point $x_1$  the inequality 
\begin{align}\label{3.31}
0\ge&~H \Delta_f \psi -2 (|\nabla \psi |^2/\psi) H + \psi \Delta_f H\nonumber\\
 \ge&-H [c_2 +(m-1) c_1(1+R \sqrt{k})+2c_1^2]/R^2 \nonumber\\
&+ \psi (e^{-h}\Sigma -\Sigma_u ) H+ \psi [2{(\Delta_f h)^2}/{m} -2 \langle \nabla h , \nabla H \rangle ]\nonumber\\
& + \psi (e^{-h}\Sigma -\Sigma_u + \mu e^h \Sigma_{uu} -2(m-1)k ) |\nabla h|^2 \nonumber\\
& -2 \psi\langle \nabla h , e^{-h} \Sigma_x - \mu \Sigma_{xu} \rangle + \psi \mu e^{-h}\Delta_f \Sigma^{x}.
\end{align}
Referring now to the above inequality, since we have $H>0$ and $\nabla (\psi H) =0$ at $x_1$, it is easily seen that,   
\begin{align} \label{3.32}
\psi \langle \nabla h, \nabla H \rangle = -H \langle \nabla h , \nabla \psi \rangle 
\le H |\nabla h| |\nabla \psi| \le c_1 \frac{\sqrt \psi}{R} H |\nabla h|. 
\end{align} 
Likewise by an application of the Cauchy-Schwarz inequality we can write 
\begin{align} \label{3.33}
\langle \nabla h , e^{-h} \Sigma_x - \mu \Sigma_{xu} \rangle \le |\nabla h| |e^{-h} \Sigma_x - \mu \Sigma_{xu}|.
\end{align}

Therefore by substituting \eqref{3.32} and \eqref{3.33} back in \eqref{3.31}, making note of \eqref{3.1} and multiplying 
through by $\psi \ge 0$, it follows that 
\begin {align}\label{3.34}
0 \ge&-\psi H [c_2 +(m-1) c_1(1+R \sqrt{k})+2c_1^2]/R^2 \nonumber\\
& + (2 \psi ^2/m) (|\nabla h|^2 + e^{-h}\Sigma)^2 - 2c_1 \psi^{3/2} |\nabla h| H/R\nonumber\\
& + \psi ^2 [ e^{-h} \Sigma - \Sigma_u + \mu e^ h \Sigma_{uu} -2(m-1)k ] |\nabla h|^2 \nonumber\\
& + \psi^2 H (e^{-h}\Sigma -\Sigma_u) - 2 \psi^2 |e^{-h} \Sigma_x - \mu \Sigma_{xu}| |\nabla h| \nonumber \\
& + \psi^2 \mu e^{-h}\Delta_f \Sigma^{x}.
\end{align}

Now in order to obtain the desired bounds out of this it is more efficient to proceed by setting 
\begin {align}  
y= \psi |\nabla h|^2, \qquad z=- \psi e^{-h}\Sigma, 
\end{align}
In particular note that $y-z=-\psi \Delta_f h$ and $y-\mu z= \psi H$ by \eqref{3.1} and \eqref{3.2} respectively. 
Substituting the above in \eqref{3.34} thus gives
\begin {align}\label{3.37}
0 \ge& -\psi H [c_2 +(m-1) c_1(1+R \sqrt{k})+2c_1^2]/R^2  \nonumber\\
&+ (2/m) [ (y-z)^2 - (mc_1/R) y^{1/2} (y-\mu z) -m\mathsf{A}_\Sigma y -m \mathsf{B}_\Sigma y^{1/2}]\nonumber\\
&+ \psi^2 H (e^{-h}\Sigma -\Sigma_u) + \psi^2 \mu e^{-h}\Delta_f \Sigma^{x}, 
\end{align}
where 
\begin{align}\label{3.38}
\mathsf{A}_\Sigma = (m-1)k - \inf_{\mathcal{B}_{2R}} \{ (e^{-h}\Sigma -\Sigma_u + \mu e^h \Sigma_{uu})_- /2 \},
\end{align}
\begin{align}\label{3.39}
\mathsf{B}_\Sigma = \sup_{\mathcal{B}_{2R}} | e^{-h} \Sigma_x -\mu \Sigma_{xu}|.
\end{align}

Utilising Lemma \ref{LiYau-basiclemma} upon setting $a=mc_1/R, b=m\mathsf{A}_\Sigma$ and 
$c= m\mathsf{B}_\Sigma$ it follows [{\it see} \eqref{eq6.9}] that  
\begin {align} \label{3.45}
(y-z)^2 - mc_1 y^{1/2}&(y-\mu z)/R - m\mathsf{A}_\Sigma y - m \mathsf{B}_\Sigma y^{1/2} \nonumber\\
\geq &\frac{1}{\mu^2}(y-\mu z)^2-\frac{m^2 c_1^2 \mu^2}{8(\mu -1)R^2}(y-\mu z)\nonumber\\
& -\frac{m^2 \mu^2 \mathsf{A}_\Sigma^2}{4(1-\varepsilon)(\mu-1)^2} 
-\frac{3}{4} \left[ \frac{m^4 \mu^2 \mathsf{B}_\Sigma^4}{4\varepsilon (\mu -1)^2}\right]^{1/3}. 
\end{align} 
Thus from \eqref{3.37} it follows that 
\begin {align}
0 \ge& -\psi H [c_2 +(m-1) c_1(1+R \sqrt{k})+2c_1^2]/R^2  \nonumber\\
&+\frac{2}{m} \bigg[ \frac{(\psi H)^2}{\mu ^2}-\frac{m^2 c_1 ^2 \mu^2}{8(\mu-1)R^2}(\psi H)\nonumber\\ 
&-\frac{m^2 \mu^2\mathsf{A}_\Sigma^2}{4(1-\varepsilon)(\mu-1)^2} 
-\frac{3}{4} \left[ \frac{m^4 \mu^2 \mathsf{B}_\Sigma^4}{4\varepsilon (\mu -1)^2}\right]^{1/3} \bigg] \nonumber\\
&+ \psi^2 H (e^{-h}\Sigma -\Sigma_u) + \psi^2 \mu e^{-h}\Delta_f \Sigma^{x}, 
\end{align}
or after basic considerations and a rearrangement of terms 
\begin {align}\label{3.46}
0 \ge&~\frac{2}{m\mu^2}(\psi H)^2 \nonumber \\
&- \bigg[ \frac{1}{R^2}[c_2 +(m-1) c_1(1+R \sqrt{k})+2c_1^2] \nonumber\\
&- \inf_{\mathcal{B}_{2R}} \{ (e^{-h}\Sigma -\Sigma_u)_{-} \} + \frac{m c_1 ^2 \mu^2}{4(\mu-1)R^2}\bigg] \psi H \nonumber \\
& -\bigg[\frac{m \mu^2\mathsf{A}_\Sigma^2}{2(1-\varepsilon)(\mu-1)^2} 
+ \frac{3}{2} \left[ \frac{m \mu^2 \mathsf{B}_\Sigma^4}{4\varepsilon (\mu -1)^2}\right]^{1/3} 
- \mu \inf_{\mathcal{B}_{2R}} \{ (e^{-h}\Delta_f \Sigma^{x})_- \} \bigg].
\end{align}
Here we have used $\psi H (e^{-h}\Sigma -\Sigma_u)_- \le \psi ^2 H [e^{-h}\Sigma -\Sigma_u]$ and 
$\mu (e^{-h}\Delta_f \Sigma^{x})_- \le \psi^2 \mu e^{-h}\Delta_f \Sigma^{x}$. Now upon setting 
\begin{align}\label{3.47}
\mathsf{D}= & ~ [c_2 +(m-1) c_1(1+R \sqrt{k})+2c_1^2]/R^2 \nonumber\\
&- \inf_{\mathcal{B}_{2R}} \{ (e^{-h}\Sigma -\Sigma_u)_{-} \} + [m c_1 ^2 \mu^2/(4(\mu-1)R^2)], 
\end{align}
and
\begin{align}\label{3.48}
\mathsf{E}=&~ m \mu^2 \mathsf{A}^2/(2(1-\varepsilon)(\mu-1)^2) \nonumber\\
&+ 3/2 [m \mu^2 \mathsf{B}^4/(4\varepsilon (\mu -1)^2)]^{1/3} 
- \mu \inf_{\mathcal{B}_{2R}} \{ (e^{-h}\Delta_f \Sigma^{x})_- \}, 
\end{align}
we can write \eqref{3.46} as 
\begin{align} \label{3.49}
0 \ge 2(\psi H)^2/(m \mu^2) -\mathsf{D} (\psi H) - \mathsf{E}.
\end{align}
As a result it follows from this inequality that 
\begin {align}
\psi H &\le (m \mu^2)/4 
\left[ \mathsf{D}+\sqrt{\mathsf{D}^2+(8 \mathsf{E})/(m\mu^2)} \right] \nonumber\\
&\le (m \mu^2)/4 \left[ 2\mathsf{D}+\sqrt{(8\mathsf{E})/(m \mu^2)} \right] 
= (m \mu^2/2) \mathsf{D} + \mu \sqrt{m\mathsf{E}/2}.
\end{align}  
Since $\psi \equiv 1$ on $\mathcal{B}_R$ and $x_1$ is a maximum point of $\psi H$ 
on $\mathcal{B}_{2R}$, we have
\begin{align}
\sup_{\mathcal{B}_R} H = \sup_{\mathcal{B}_R} [\psi H] \le \sup_{\mathcal{B}_{2R}} [\psi H] = (\psi H)(x_1).
\end{align}
Thus it follows that  
\begin{align}
\sup_{\mathcal{B}_R} H \le (m \mu^2/2)\mathsf{D} + \mu \sqrt{m \mathsf{E}/2}.
\end{align} 
Therefore recalling \eqref{3.2}, substituting for $\mathsf{D}$ and $\mathsf{E}$ from \eqref{3.47} and \eqref{3.48} 
above, we can write after multiplying both sides by $1/\mu$, that for every $x \in \mathcal{B}_R$: 
\begin{align}
\mu^{-1} |\nabla h|^2 + e^{-h} \Sigma (x,e^h) \le& ~m \mu \mathsf{D}/2 + (m\mathsf{E}/2)^{1/2} \nonumber\\
\le& ~m \mu [c_2 +(m-1) c_1(1+R \sqrt{k})+2c_1^2]/(2R^2) \nonumber\\
&- (m\mu/2) \inf_{\mathcal{B}_{2R}} \{ (e^{-h}\Sigma -\Sigma_u)_{-} \} \nonumber\\
&+ [m^2 c_1 ^2 \mu^3/(8(\mu-1)R^2)]\nonumber\\
& + \sqrt{m} \{[m \mu^2 \mathsf{A}_\Sigma^2/(4(1-\varepsilon)(\mu-1)^2)] \nonumber\\
&+ (3 /4) [m \mu^2 \mathsf{B}_\Sigma^4/(4\varepsilon (\mu -1)^2)]^{1/3} \nonumber\\
&- (\mu/ 2 )\inf_{\mathcal{B}_{2R}} \{ (e^{-h}\Delta_f \Sigma^{x})_- \}\}^{1/2}.
\end{align}
Finally reverting back to $u$ upon noting the relation $h = \log u$ and rearranging terms
\begin{align}
\frac{|\nabla u|^2}{\mu u^2} + \frac{\Sigma [x, u]}{u} \le& ~
m \mu [c_2 +(m-1) c_1(1+R \sqrt{k})\nonumber\\
&+2c_1^2 + mc_1^2 \mu^2/(4(\mu-1) )]/ (2R^2) \nonumber\\
&- (m\mu/2) \inf_{\mathcal{B}_{2R}} \{ (\Sigma/u -\Sigma_u)_{-} \} \nonumber\\
& + \sqrt{m} \{[m \mu^2 \mathsf{A}_\Sigma^2/(4(1-\epsilon)(\mu-1)^2)] \nonumber\\
&+ (3 /4) [m \mu^2 \mathsf{B}_\Sigma^4/(4\epsilon (\mu -1)^2)]^{1/3} \nonumber\\
&- (\mu/ 2 )\inf_{\mathcal{B}_{2R}} \{ ([\Delta_f \Sigma^{x}]/u)_- \}\}^{1/2},
\end{align}
which is the desired estimate as in \eqref{1.26}. The proof is thus complete.

\section{Proof of the elliptic Harnack inequality in Theorem \ref{cor Harnack}} 
\label{Harnack-Section}

Now in order to prove the Harnack inequality we need to integrate the differential Harnack inequality 
along a geodesic path $\gamma$ joining the points $x_1$ and $x_2$ inside $\mathcal{B}_{R}$.
Towards this end let us begin by rewriting the local gradient estimate for \eqref{eq11} as follows: 

\begin{align}\label {4.1}
\sup_{\mathcal{B}_R} \frac{|\nabla u|^2}{u^2}  \le&~\frac{m\mu^2}{2R^2}\left[ [c_2 +(m-1) c_1(1+R \sqrt{k})+2c_1^2] 
+ \frac{mc_1^2 \mu^2}{4(\mu-1)}\right]\nonumber\\
& +\frac{\sqrt {m} \mu}{2}\bigg[ \frac{m \mu^2\mathsf{A}^2}{(1-\varepsilon)(\mu-1)^2} 
+ \left[\frac{27m \mu^2\mathsf{B}^4}{ 4\varepsilon(\mu-1)^2}\right]^{1/3} 
- 2\mu \inf_{\mathcal{B}_{2R}} \left\{ \frac{(\Delta_f \Sigma^x)_-}{u} \right\} \bigg] ^{1/2} \nonumber \\
&+ \frac{m \mu^2}{2} \sup_{\mathcal{B}_{2R}} \left\{ \frac{(u \Sigma_u[x,u]-\Sigma [x,u])_+}{u} \right\} 
- \mu \inf_{\mathcal{B}_{R}} \left\{\frac{(\Sigma [x,u])_-}{u} \right\} := \mathbb{H}.
\end{align}

Here we have denoted the expression on the right-hand side of \eqref{4.1} by $\mathbb{H}$ which is a positive 
constant. Now integrating the quantity $ |\nabla u|/u$ along a geodesic curve $\gamma$ in $\mathcal{B}_R$ 
(with $\gamma (0)=x_1$ and $\gamma(1) = x_2$) we have
\begin{align}
\log u(x_2) - \log u(x_1)&= \int_0^1 \frac{d}{ds} \log u(\gamma(s)) \, ds \nonumber \\
&= \int_0^1 \langle [\nabla u/u](\gamma(s)) , \gamma'(s) \rangle \, ds 
\le \left[ \sup_{\mathcal{B}_{R}}  \frac{|\nabla u|}{u} \right] \int_0^1 |\gamma'| \, ds \nonumber\\
& \le d(x_1, x_2) \sqrt \mathbb{H}  \le 2R \sqrt \mathbb{H}.
\end{align}
Therefore $\log [u(x_2)/u(x_1)] \le d(x_1, x_2) \sqrt \mathbb{H} \le 2R \sqrt \mathbb{H}$ or after exponentiation 
\begin{equation}
u(x_2) \le e^{2R \sqrt \mathbb{H} } u(x_1)
\end{equation}
giving the desired inequality. The remaining assertions are now straightforward consequences of this inequality. \hfill $\square$

\section{Proof of the Liouville result in Theorem \ref{thm Liouville}} 
\label{Liouville-Section}

Starting from \eqref{1.26-global} and noting that $\mathsf{B}_\Sigma =0$ (as a result of $|\Sigma_x -\mu u \Sigma_{xu}| \equiv 0$), 
$k=0$ and $\mathsf{C}_\Sigma = 0$ (as a result of $\Delta_f \Sigma^x \equiv 0$) we obtain, after rearranging the inequality,  
\begin{align} \label{eq2.9}
\frac{|\nabla u|^2}{\mu u^2} + \frac{\Sigma(u)}{u} \le&~ 
\frac{\sqrt{m}}{2} \bigg[ \frac{m \mu^2\mathsf{A}_\Sigma^2}{(1-\varepsilon)(\mu-1)^2} 
+ \left[\frac{27 m \mu^2\mathsf{B}_\Sigma^4}{4 \varepsilon(\mu-1)^2}\right]^{1/3}
+ 2 \mu \mathsf{C}_\Sigma \bigg]^{1/2} \nonumber\\ 
& + \frac{m \mu}{2} \sup_{M}\left\{ \frac{(u\Sigma_u [u]-\Sigma[u])_+}{u} \right \} \nonumber \\
\le&~ m \mu \bigg[ 
\sup_{M} \left\{ \frac{(-\Sigma[u] + u\Sigma_u[u] - \mu u^2 \Sigma_{uu}[u])_+}{4(\sqrt{1-\varepsilon})(\mu-1) u} \right\} \nonumber \\
& + \sup_{M}\left\{ \frac{(u\Sigma_u[u]-\Sigma[u])_+}{2u} \right \} \bigg]. 
\end{align}
Next from the imposed assumptions on $\Sigma$ and its derivatives it is easily seen that 
\begin{equation}
(-\Sigma[u] + u\Sigma_u[u] - \mu u^2 \Sigma_{uu}[u])_+ \equiv 0, \qquad (u\Sigma_u[u]-\Sigma[u])_+ \equiv 0.
\end{equation}
Hence from \eqref{eq2.9} it follows that 
\begin{equation}
\frac{|\nabla u|^2}{\mu u^2} + \frac{\Sigma[u]}{u} \equiv 0,
\end{equation}
and so again from the assumptions imposed on $\Sigma$ that $|\nabla u|^2/u^2 \equiv 0$. The conclusion now follows at once. \hfill $\square$

\qquad \\\
{\bf Authors Statements.} The authors declare no conflict of interest. They have equal contribution in this research. The datasets generated 
during and/or analysed during the current study are available from the corresponding author on reasonable request. The authors gratefully 
acknowledge support from EPSRC.

\end{document}